# On the Finite Number of Directional Stationary Values of Piecewise Programs

Ying Cui [*]    Jong-Shi Pang [†]

**Abstract**

Extending a fundamental result for (indefinite) quadratic programs, this paper shows that certain non-convex piecewise programs have only a finite number of directional stationary values, and thus, possess only finitely many locally minimum values. We present various special cases of our main results, in particular, an application to a least-squares piecewise affine regression problem for which every directional stationary point is locally minimizing.

## 1 Introduction.

Our study is motivated by some empirical results that we obtained in computational experiments for solving the non-traditional least-squares piecewise affine regression problems formulated as non-convex, non-differentiable optimization problems of the piecewise quadratic type [8]. Applying a majorization-minimization based algorithm to these problems and starting at hundreds of initial iterates, we observed that the final objective values at termination of the algorithm clustered around a few distinguished values. By the convergence theory of the algorithm, we know that such values are the objective function values at (approximate) *directional stationary solutions* of the problems, which are stationary solutions defined by the elementary directional derivatives of these piecewise functions. This computational observation is reminiscent of a classic result proved by Luo and Tseng in [22, Lemma 3.1] for standard (indefinite) quadratic programs and begs the question of whether a linearly constrained "piecewise quadratic program" indeed has only finitely many *directional stationary values*. It turns out that this question has easy answers for certain classes of such problems, not-so-easy answers for others, and difficult and indeed un-resolved answers in the general sense of the question. The goal of this paper is to properly address the question, provide partial answers to it, and leave open its full resolution as a conjecture. The results obtained here may serve as a first step to analyze the landscape of non-convex piecewise programs, and have relevance to the sequential convergence of objective values in iterative descent algorithms for solving these non-differentiable optimization problems. We will illustrate these points with the class of least-squares piecewise affine regression problems; see the last section of this paper, in particular, Proposition 14.

[*]The Daniel J. Epstein Department of Industrial and Systems Engineering, University of Southern California, Los Angeles, California, USA. Email: yingcui@usc.edu

[†]The Daniel J. Epstein Department of Industrial and Systems Engineering, University of Southern California, Los Angeles, California, USA. Email: jongship@usc.edu

# 2  Piecewise Programs and Directional Stationarity.

In this section, we summarize some preliminary materials needed for the rest of the paper. A continuous function $\psi : \mathbb{R}^n \to \mathbb{R}$ is said to be *piecewise affine* (PA) if there exist finitely many affine functions $\{\psi_1, \ldots, \psi_k\}$ for some positive integer $k$ such that $\psi(x) \in \{\psi_1(x), \ldots, \psi_k(x)\}$ for all $x \in \mathbb{R}^n$. It is known from [36] that every piecewise affine function has a max-min representation which can be equivalently expressed in following difference-max form:

$$\psi(x) = \max_{1 \leq i \leq k_1} \left( (a^i)^T x + \alpha_i \right) - \max_{1 \leq j \leq k_2} \left( (b^j)^T x + \beta_j \right) \tag{1}$$

for some $n$-dimensional vectors $\{a^i\}_{i=1}^{k_1}$ and $\{b^j\}_{j=1}^{k_2}$ and scalars $\{\alpha_i\}_{i=1}^{k_1}$ and $\{\beta_j\}_{j=1}^{k_2}$. This difference-max representation provides the basis for extensions to a broad class of piecewise quadratic (PQ) functions to be studied in this paper.

A continuous function $\psi : \Omega \subseteq \mathbb{R}^n \to \mathbb{R}$ is said to be *piecewise quadratic* (PQ) if there exist finitely many quadratic functions $\{\psi_1, \ldots, \psi_k\}$ for some positive integer $k$ such that $\psi(x) \in \{\psi_1(x), \ldots, \psi_k(x)\}$ for all $x \in \Omega$. Clearly, a function of the form

$$\psi(x) = \max_{1 \leq i \leq k_1} \psi_{1i}(x) - \max_{1 \leq i \leq k_2} \psi_{2i}(x), \tag{2}$$

where each $\psi_{ji}$ is a quadratic function, is PQ; nevertheless, unlike the PA functions, it is not known if every PQ function has the above representation as the difference-max of quadratic functions. Although the Ph.D. thesis [40] and the subsequent paper [41] have studied, respectively, piecewise quadratic programming and structures of convex piecewise quadratic functions extensively, this issue has not been addressed fully in the literature to date. Setting aside this question or representation, we will devote much of our analysis to the class of PQ functions with the representation (2). A subclass of PQ functions consists of the *piecewise linear-quadratic* (PLQ) functions [34, Definition 10.20] whose domain $\Omega$ can be represented as the union of finitely many polyhedral sets on each of which $\psi(x)$ is a quadratic function. As noted in the latter reference with the piecewise quadratic function in 2 variables: $\psi(x_1, x_2) = |x_1^2 + x_2^2 - 1|$, a PQ function may not be a PLQ function. One may refer to [40, 41] for more properties of the PLQ functions. Incidentally, it is known that if $\psi$ is once continuously differentiable on an open set, then it is piecewise quadratic there if and only if it is piecewise linear-quadratic; see [7, Proposition 3.1].

Let the domain $\Omega$ be an open subset of $\mathbb{R}^n$. The directional derivative of $\psi$ at a vector $\bar{x} \in \Omega$ in the direction $v \in \mathbb{R}^n$ is given by

$$\psi'(\bar{x}; v) \triangleq \lim_{\delta \downarrow 0} \frac{\psi(\bar{x} + \delta v) - \psi(\bar{x})}{\delta},$$

if it exists. We say that $\psi$ is *directionally differentiable* (dd) on $\Omega$ if the directional derivatives $\psi'(x; v)$ exist for all $(x, v) \in \Omega \times \mathbb{R}^n$. Let $X$ be a closed convex set in $\Omega$ and $\psi$ be a dd function. We say that $\bar{x} \in X$ is a *d(irectional)-stationarity point* of the program

$$\underset{x \in X}{\text{minimize}} \; \psi(x), \tag{3}$$

if $\psi'(\bar{x}, x - \bar{x}) \geq 0$ for all $x \in X$, or equivalently,

$$0 \in \widehat{\partial}\psi(\bar{x}) + \mathcal{N}(\bar{x}; X),$$

where

$$\widehat{\partial}\psi(\bar{x}) \triangleq \left\{ v \in \mathbb{R}^n \;\middle|\; \liminf_{\substack{x \to \bar{x} \\ x \neq \bar{x}}} \frac{\psi(x) - \psi(\bar{x}) - v^T(x - \bar{x})}{\| x - \bar{x} \|} \geq 0 \right\}$$

denotes the *regular subdifferential* of $\psi$ at $\bar{x}$ [34, Definition 8.3] and $\mathcal{N}(\bar{x}; X)$ denotes the normal cone of the convex set $X$ at $\bar{x}$ as in convex analysis. For the pair $(\psi, X)$ as given, d-stationarity of $\bar{x}$ is necessary for $\bar{x}$ to be a local minimizer of (3) [34, Theorem 10.1]. Throughout this paper, we let $\mathcal{D}_{(\psi, X)}$ denote the set of all d-stationary points of (3). For any $x \in \mathcal{D}_{(\psi, X)}$, we call $\psi(x)$ a *d-stationary value* of $\psi$ on $X$. The basic result of Luo and Tseng [22, Lemma 3.1] states that $\psi\left(\mathcal{D}_{(\psi, X)}\right)$ is a (possibly empty) finite set if $\psi$ is a quadratic function and $X$ is a polyhedron. This is the result that we aim to extend to the case where $\psi$ is piecewise quadratic. It should be mentioned that finiteness does not imply polynomiality; in fact, the number of such values can be expected to be exponential in general; this is already clear in the Luo-Tseng result for quadratic programs and will become more evident in the proof of the new results.

For any real-valued function $\psi$, $\widehat{\partial}\psi(\bar{x})$ is a subset of $\partial_C \psi(\bar{x})$ [24, Theorem 3.57], where the latter is the Clarke subdifferential of $\psi$ at $\bar{x}$ that defined as

$$\partial_C \psi(\bar{x}) \triangleq \left\{ v \in \mathbb{R}^n \ \Big| \ \limsup_{\substack{x \to \bar{x} \\ t \downarrow 0}} \frac{\psi(x + tw) - \psi(x) - t\, v^T w}{t} \geq 0, \quad \forall\, w \in \mathbb{R}^n \right\}.$$

Consider a dd function $\psi$ given by

$$\psi(x) = \psi_1(x) - \max_{1 \leq i \leq k_2} \psi_{2i}(x), \tag{4}$$

where $\psi_1$ and each $\psi_{2i}$ are dd functions defined on the same domain $\Omega$. Since for any $\bar{x} \in \Omega$,

$$\psi'(\bar{x}; d) = \psi_1'(\bar{x}; d) - \max_{i \in \mathcal{A}_2(\bar{x})} \psi_{2i}'(\bar{x}; d), \quad \text{where } \mathcal{A}_2(\bar{x}) \triangleq \operatornamewithlimits{argmax}_{1 \leq i \leq k_2} \psi_{2i}(\bar{x}),$$

we easily obtain the following lemma.

**Lemma 1.** Let $\psi$ be given by (4). Then $\bar{x}$ is a d-stationary point of (3) if and only if $\bar{x}$ is a d-stationary point of $\operatornamewithlimits{minimize}_{x \in X} \psi_1(x) - \psi_{2; \bar{i}}(x)$ for any $\bar{i} \in \mathcal{A}_2(\bar{x})$.

The important point about this lemma is that the non-differentiable pointwise maximum of the $\psi_{2i}$ functions has been replaced by the individual maximands in d-stationarity; this equivalence is particularly useful if each $\psi_{2i}$ is differentiable, such as quadratic, as we will see in the later development. This result highlights an essential property of directional stationarity, which may not be possessed by other relaxed stationarity concepts such as that based on the Clarke subdifferential.

## 3 Linearly Constrained Piecewise Programs.

In this section, we focus on the case where the constraint set $X$ in (3) is a polyhedron and show that the problem (3) has finitely many d-stationary values, and thus finitely many local minimum values, for several classes of piecewise functions $\psi$. The following result is not difficult to prove. Each of the three types of functions in this result is different from the other two.

**Proposition 2.** Suppose that $X$ is a polyhedral set. The set $\psi\left(\mathcal{D}_{(\psi, X)}\right)$ is finite if any one of the following three conditions hold:

(a) $\psi$ is a piecewise linear-quadratic function on $X$;

(b) $\psi$ is given by (4) where $\psi_1$ is convex and each $\psi_{2i}$ is concave for all $i = 1, \cdots, k_2$;

(c) $\psi(x) = \psi_1(x) + \psi_2(x)$ with $\psi_1$ being a convex function and $\psi_2$ being a piecewise affine function.

**Proof.** (a) Let $\psi$ equal to a quadratic function $q_i(x)$ on the polyhedral set $C_i$ for $i = 1, \cdots, k$ whose union is $X$. For any $\bar{x} \in X$, let $\mathcal{A}(\bar{x}) \subseteq \{1, \ldots, k\}$ be the set of indices $i$ such that $\bar{x} \in C_i$. It is known from (the proof of) [34, Proposition 10.21] that

$$\psi'(\bar{x}; x - \bar{x}) = \nabla q_i(\bar{x})^T (x - \bar{x}) \quad \text{if } x \in C_i, \, i \in \mathcal{A}(\bar{x}).$$

For $i = 1, \ldots, k$, let $\mathcal{D}_i$ denote the set of all d-stationary points of $q_i(x)$ on $C_i$; i.e.,

$$\mathcal{D}_i = \left\{ \bar{x} \in C_i \mid \nabla q_i(\bar{x})^T (x - \bar{x}) \geq 0, \, \forall \, x \in C_i \right\}.$$

Thus, $\bar{x} \in \mathcal{D}_{(\psi, X)}$ implies that $\bar{x} \in \bigcap_{i \in \mathcal{A}(\bar{x})} \mathcal{D}_i$. Therefore, $\mathcal{D}_{(\psi, X)} \subseteq \bigcup_{i=1}^{k} \mathcal{D}_i$. By [22, Lemma 3.1], each $q_i(\mathcal{D}_i)$ is a finite set. Since $\psi\left(\mathcal{D}_{(\psi, X)}\right) \subseteq \bigcup_{i=1}^{k} q_i(\mathcal{D}_i)$, the finiteness of $\psi\left(\mathcal{D}_{(\psi, X)}\right)$ follows readily.

To prove the same finiteness under condition (b), we note that by Lemma 1, $\bar{x}$ is a d-stationary point of $\psi$ on $X$ if and only if $\bar{x}$ is a d-stationary point of the convex function $\varphi_i(x) \triangleq \psi_1(x) - \psi_{2i}(x)$ on $X$ for every $i \in \mathcal{A}_2(\bar{x})$. Since the latter index set is finite, it follows that the set $\psi(\mathcal{D}_{(\psi, X)})$ is a subset of the finite set of global minimum values $\left\{ \underset{x \in X}{\text{minimum}} \, \varphi_i(x) \right\}_{i=1}^{k_2}$.

Finally, statement (c) is a special case of (b) because by the difference-max representation (1), we have $\psi_2(x) = \widehat{\psi}_2(x) - \max_{1 \leq i \leq k_2} \psi_{2i}(x)$, where $\widehat{\psi}_2$, being the pointwise maximum of finitely many affine functions, is a convex function, and each $\psi_{2i}$ is an affine, thus concave function. $\square$

**Remark 3.** *A weaker notion of stationarity is called criticality in the difference-of-convex (dc) literature [30]. Suppose that the function $\psi$ in the program (3) is the difference of two convex functions $\psi(x) = \psi_1(x) - \psi_2(x)$. Then $x \in X$ is called a critical point of $\psi$ on $X$ if*

$$(\partial \psi_1(x) + \mathcal{N}(x; X)) \cap \partial \psi_2(x) \neq \emptyset.$$

*It is possible for a dc function in the form of the above proposition to have infinitely many critical values (i.e., the objective values at the critical points). Consider a univariate convex function $\psi(x) = \max(x, 0)$ for $x \in \mathbb{R}$. Obviously this function has a unique d-stationary value $= 0$, which is the globally minimum value of $\psi$; nevertheless every $x \leq 0$ is a global minimizer, thus d-stationary point. This simple function can be written as any of the forms (a)–(c) in Proposition 2. Under the dc decomposition $\psi(x) = \psi_1(x) - \psi_2(x)$, where*

$$\psi_1(x) = \begin{cases} (2n+1)x - n(n+1) & \text{if } x \in [n, n+1), \quad n = 0, 1, 2, \cdots, \\ 0 & \text{if } x \leq 0 \end{cases}$$

*and*

$$\psi_2(x) = \begin{cases} 2nx - n(n+1) & \text{if } x \in [n, n+1), \quad n = 0, 1, 2, \cdots, \\ 0 & \text{if } x \leq 0, \end{cases}$$

*one may verify that*

$$\partial \psi_1(n) = [2n-1, 2n+1] \quad \text{and} \quad \partial \psi_2(n) = [2n-2, 2n].$$

*Therefore,*

$$\partial \psi_1(n) \cap \partial \psi_2(n) = [2n-1, 2n] \neq \emptyset.$$

This shows that $x = n$ is a critical point of $\psi_1(x) - \psi_2(x)$ for any nonnegative integer $n$, leading to infinitely many critical values $\psi(n) = n$. See Figure 1 for an illustration of the function $\psi$ and its dc representation.

The above example confirms that the critical point, which depends on the particular dc representation of the objective function, is a weaker concept than the d-stationary point. More seriously, for a convex function, such as the simple univariate function here, a critical point may not be a global minimizer if the dc decomposition is chosen improperly. □

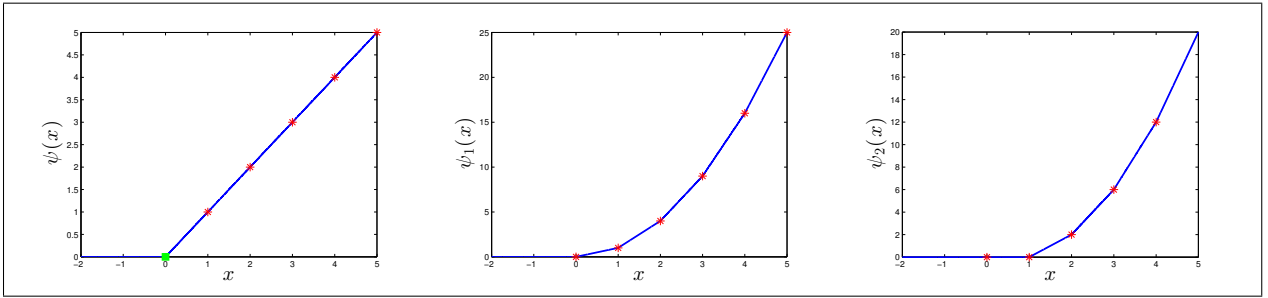

Figure 1: the convex function $\psi$ and its dc representation

The next result pertains to a difference-max function of the form (2) where each $\psi_{ji}$ is a quadratic function. It identifies subsets $S$ of $X$ for which the set $\psi\left(S \cap \mathcal{D}_{(\psi,X)}\right)$ is finite.

**Proposition 4.** Suppose that $X$ is a polyhedral set. Let $\psi$ be given by (2) where each $\psi_{ji}$ is a quadratic function. For an arbitrary tuple $c \triangleq \left(c^i\right)_{i=1}^{k_1}$ of $n$-vectors, let

$$S_c \triangleq \left\{ x \in X \mid \nabla \psi_{1i}(x) + c^i \text{ is a common vector for all } i \in \mathcal{A}_1(x) \right\},$$

where $\mathcal{A}_1(x) \triangleq \underset{1 \le i \le k_1}{\operatorname{argmax}} \psi_{1i}(x)$. It holds that the set $\psi\left(S_c \cap \mathcal{D}_{(\psi,X)}\right)$ is finite.

**Proof.** By Lemma 1, it suffices to prove the proposition for $\psi(x) = \underset{1 \le \ell \le k_1}{\max} \psi_{1\ell}(x) - q_2(x)$ where $q_2$ and each $\psi_{1\ell}$ are quadratic functions, the former with a constant Hessian matrix $Q$ and the latter with a constant Hessian matrix $P^\ell$. Write the polyhedral set $X \triangleq \{x \in \mathbb{R}^n \mid Ax \le b\}$ for some matrix $A \in \mathbb{R}^{m \times n}$ and vector $b \in \mathbb{R}^m$. Let $\bar{x}$ be a d-stationary point of $\psi$ on $X$. It then follows there must exist nonnegative scalars $\{\lambda_\ell^{\bar{x}}\}_{\ell \in \mathcal{A}_1(\bar{x})}$ summing to unity and nonnegative scalars $\{\mu_j^{\bar{x}}\}_{j \in \mathcal{A}_X(\bar{x})}$, where $\mathcal{A}_X(\bar{x}) \triangleq \{j \mid A_{j\bullet}\bar{x} = b_j\}$ with $A_{j\bullet}$ denoting the $j$-th row of the matrix $A$ is the index set of active constraints at $\bar{x}$, such that

$$\sum_{\ell \in \mathcal{A}_1(\bar{x})} \lambda_\ell^{\bar{x}} \nabla \psi_{1\ell}(\bar{x}) - \nabla q_2(\bar{x}) + \sum_{j \in \mathcal{A}_X(\bar{x})} \mu_j^{\bar{x}} \left( A_{j\bullet} \right)^T = 0.$$

Note that both sets of multipliers: $\{\lambda_\ell^{\bar{x}}\}_{\ell \in \mathcal{A}_1(\bar{x})}$ and $\{\mu_j^{\bar{x}}\}_{j \in \mathcal{A}_X(\bar{x})}$ are dependent on the d-stationary point $\bar{x}$. Let $y$ be another d-stationary point of $\psi$ on $X$ such that $\mathcal{A}_1(y) = \mathcal{A}_1(\bar{x}) = \mathcal{I}$ and $\mathcal{A}_X(y) = \mathcal{A}_X(\bar{x}) = \mathcal{J}$ for some common index sets $\mathcal{I}$ and $\mathcal{J}$, respectively. We evaluate

$$\begin{aligned}
\psi(y) - \psi(\bar{x}) &= \sum_{\ell \in \mathcal{I}} \lambda_\ell^{\bar{x}} \left( \psi_{1\ell}(y) - \psi_{1\ell}(\bar{x}) \right) - \left( q_2(y) - q_2(\bar{x}) \right) \\
&= \sum_{\ell \in \mathcal{I}} \lambda_\ell^{\bar{x}} \left[ \nabla \psi_{1\ell}(\bar{x})^T (y - \bar{x}) + \tfrac{1}{2} (y - \bar{x})^T P^\ell (y - \bar{x}) \right] - \\
&\qquad \left[ \nabla q_2(\bar{x})^T (y - \bar{x}) + \tfrac{1}{2} (y - \bar{x})^T Q (y - \bar{x}) \right] \\
&= -\sum_{j \in \mathcal{J}} \mu_j^{\bar{x}} A_{j\bullet} (y - \bar{x}) + \tfrac{1}{2} \sum_{\ell \in \mathcal{I}} \lambda_\ell^{\bar{x}} (y - \bar{x})^T P^\ell (y - \bar{x}) - \tfrac{1}{2} (y - \bar{x})^T Q (y - \bar{x}) \\
&= \tfrac{1}{2} \sum_{\ell \in \mathcal{I}} \lambda_\ell^{\bar{x}} (y - \bar{x})^T P^\ell (y - \bar{x}) - \tfrac{1}{2} (y - \bar{x})^T Q (y - \bar{x}), \quad \text{because } A_{\mathcal{J}\bullet} y = b_{\mathcal{J}} = A_{\mathcal{J}\bullet} x.
\end{aligned}$$

If $\bar{x}$ and $y$ are both in $S_c$, then for any two indices $\ell$ and $\ell'$ in $\mathcal{I}$, we have

$$\nabla \psi_{1\ell}(\bar{x}) - \nabla \psi_{1\ell'}(\bar{x}) = c^{\ell'} - c^\ell = \nabla \psi_{1\ell}(y) - \nabla \psi_{1\ell'}(y);$$

which implies

$$P^\ell (y - \bar{x}) = \nabla \psi_{1\ell}(y) - \nabla \psi_{1\ell}(\bar{x}) = \nabla \psi_{1\ell'}(y) - \nabla \psi_{1\ell'}(\bar{x}) = P^{\ell'} (y - \bar{x}).$$

Thus, for any $\bar{\ell} \in \mathcal{I}$, we deduce

$$\psi(y) - \psi(\bar{x}) = \tfrac{1}{2} (y - \bar{x})^T \left[ P^{\bar{\ell}} - Q \right] (y - \bar{x}).$$

By symmetry, we also have

$$\psi(\bar{x}) - \psi(y) = \tfrac{1}{2} (y - \bar{x})^T \left[ P^{\bar{\ell}} - Q \right] (y - \bar{x}).$$

Thus $\psi(\bar{x}) = \psi(y)$. Hence we have proved that $\psi$ is a constant on the set $S_c \cap \mathcal{D}_{(\psi, X)} \cap \mathcal{D}_{\mathcal{I}\mathcal{J}}$, where

$$\mathcal{D}_{\mathcal{I}\mathcal{J}} \triangleq \{ x \in X \mid \mathcal{A}_1(x) = \mathcal{I} \text{ and } \mathcal{A}_X(x) = \mathcal{J} \}.$$

Since there are only finitely many subsets $\mathcal{I}$ of $\{1, \cdots, k_1\}$ and $\mathcal{J}$ of $\{1, \cdots, m\}$, it follows that the set $\psi\left( S_c \cap \mathcal{D}_{(\psi, X)} \right)$ is finite. $\square$

The corollary below follows readily from Proposition 4.

**Corollary 5.** Suppose that $X$ is a polyhedral set. The following two statements hold.

(a) The set $\psi\left( \mathcal{D}_{(\psi, X)} \right)$ is a finite set if $\psi(x) = \max_{1 \leq i \leq k_1} \psi_{1i}(x) - \max_{1 \leq j \leq k_2} \psi_{2j}(x)$ with each $\psi_{1i}$ being an affine function and each $\psi_{2j}$ being an (indefinite) quadratic function.

(b) The set $\psi(\mathcal{D}_{(\psi, X)} \cap \mathcal{D}_F)$ is a finite set if $\psi(x) = \max_{1 \leq i \leq k_1} \psi_{1i}(x) - \max_{1 \leq j \leq k_2} \psi_{2j}(x)$ with $\psi_{1i}$ and $\psi_{2j}$ all being (indefinite) quadratic functions, where $\mathcal{D}_F$ is the set of all F(réchet) differentiable points of $\max_{1 \leq j \leq k_1} \psi_{1i}(x)$.

**Proof.** (a) Suppose that $\psi_{1i}(x) = (a^i)^T x + \alpha_i$ for some vectors $a^i \in \mathbb{R}^n$ and scalars $\alpha_i \in \mathbb{R}$. It follows that $\nabla \psi_{1i}(x) - a^i = 0$ is independent of $i$ for all $x \in X$.

(b) It is known that $\bar{x}$ is a F-differentiable point of $\max_{1 \leq j \leq k_1} \psi_{1i}(x)$ if and only if $\nabla \psi_{1i_1}(x) = \nabla \psi_{1i_2}(x)$ for all $i_1$ and $i_2$ in $\operatorname*{argmax}_{1 \leq i \leq k_1} \psi_{1i}(x)$ [33, Theorem 1]; equivalently, if and only if $\nabla \psi_{1i}(x)$ is independent of $i \in \operatorname*{argmax}_{1 \leq i \leq k_1} \psi_{1i}(x)$. $\square$

In Corollary 5(b), we show that the difference-max of quadratic functions has finitely many values over the set of d-stationary points at which the first max term is F-differentiable. We next move beyond such differentiable points and consider the points where there are no more than two active pieces in $\max_{1 \leq i \leq k_1} \psi_{1i}$. To prove this result, we first provide a simple determinantal lemma.

**Lemma 6.** For any two square matrices $A$ and $B$ in $\mathbb{R}^{n \times n}$ and a scalar $t \in \mathbb{R}$, the inverse of $A + tB$, if exists, can be represented as

$$(A + tB)^{-1} = \frac{1}{\det(A + tB)} \left( \sum_{k=1}^{n-1} t^k C^k + C^0 \right),$$

where $\{C^k\}_{k=0}^{n-1}$ are $n \times n$ square matrices dependent on $A$ and $B$ but not on $t$.

**Proof.** The characteristic polynomial of $A + tB$ has the form

$$p(\lambda) \triangleq \det(\lambda I_n - (A + tB)) = \lambda^n + \sum_{i=1}^{n-1} \left( \sum_{j=1}^{n-i} c_{ij} t^j + c_{i0} \right) \lambda^i + (-1)^n \det(A + tB),$$

for some scalars $c_{ij}$ depending on $A$ and $B$ only, where $I_n$ is the identity matrix of order $n$. By the Cayley-Hamilton theorem,

$$0_{n \times n} = p(A + tB) = (A + tB)^n + \sum_{i=1}^{n-1} \left( \sum_{j=1}^{n-i} c_{ij} t^j + c_{i0} \right) (A + tB)^i + (-1)^n \det(A + tB) I_n.$$

We then obtain the stated results by multiplying $(A + tB)^{-1}$ to both sides, expanding the powers $(A + tB)^i$, and rearranging terms. □

We next prove a lemma that allows us to reduce the linearly constrained problem (3) to finitely many unconstrained problems. This reduction will be used in the proof of Proposition 8 and subsequently. In essence, Luo-Tseng employed the same reduction in the proof of the finiteness result of d-stationary values of a quadratic program.

**Lemma 7.** Let $X$ be a polyhedral set in $\mathbb{R}^n$ and $\psi$ be a directionally differentiable function. If $\bar{x}$ is a d-stationary point of $\psi$ on $X$, then $\bar{x}$ is a d-stationary point of $\psi$ on the affine subspace $\{x \mid A_{I\bullet} x = b_I\}$, where $I \triangleq \mathcal{A}_X(\bar{x})$ is the index set of active constraints of $X$ at $\bar{x}$.

**Proof.** This holds because if $x$ is an element of the affine subspace in question, then $\bar{x} + \delta(x - \bar{x})$ is an element of $X$ for all $\delta > 0$ sufficiently small. □

The next result extends Corollary 5(b). The proof is quite different from that of the corollary.

**Proposition 8.** Suppose that $X$ is a polyhedral set. Let $\psi(x) = \max_{1 \leq i \leq k_1} \psi_{1i}(x) - \max_{1 \leq j \leq k_2} \psi_{2j}(x)$, where $\psi_{1i}$ and $\psi_{2j}$ are all quadratic functions on $\mathbb{R}^n$. Then the problem (3) has only finitely many d-stationary values on

$$\widetilde{\mathcal{D}} \triangleq \left\{ \bar{x} \in \mathcal{D}_{(\psi, X)} \;\middle|\; \left| \operatorname*{argmax}_{1 \leq i \leq k_1} \psi_{1i}(\bar{x}) \right| \leq 2 \right\}.$$

**Proof.** By Lemma 1 and the restriction to the set $\widetilde{\mathcal{D}}$, to prove this proposition, it suffices to show that the problem of the form

$$\operatorname*{minimize}_{x \in X} \quad \max(\psi_{1i_1}(x), \psi_{1i_2}(x)) - \psi_{2j}(x), \quad \text{where } 1 \leq i_1, i_2 \leq k_1 \text{ and } 1 \leq j \leq k_2$$

has only finitely many d-stationary values. We make two more simplifications in the following proof, without loss of generality.

(a) Write $X \triangleq \{x \in \mathbb{R}^n \mid Ax \leq b\}$ for some matrix $A \in \mathbb{R}^{m \times n}$ and vector $b \in \mathbb{R}^m$. By Lemma 7, it follows that if $\bar{x}$ is a d-stationary point of $\psi$ on $X$, then $\bar{x}$ is a d-stationary point of $\psi$ on the affine subspace $\{x \mid A_{I\bullet}x = b_I\}$, where $I \triangleq \mathcal{A}_X(\bar{x})$ is the index set of active constraints of $X$ at $\bar{x}$. As in the proof of Proposition 4, we can consider d-stationary points of $\psi$ on $X$ that have the same index sets of active constraints. By focusing on any such affine subspace with the same index set $I$, we may represent any feasible solution by $x = x^0 + Zy$, where $Z$ is a basis of the null space of $A_{I\bullet}$; we can then eliminate this linear constraint and reformulate the problem (on that subspace) in terms of the free variable $y$; see, e.g., [12, Chapter 10]. Thus, working with active constraints and employing this substitution of variables, we may take the set $X$ to be the entire space $\mathbb{R}^n$.

(b) In addition, the function $\psi_{1i_2}$ can be taken to be 0 by redefining $\psi_{1i_1}$ as $\psi_{1i_1} - \psi_{1i_2}$ and $\psi_{2j}$ as $\psi_{2j} - \psi_{1i_2}$. Therefore, it suffices to show that the unconstrained program

$$\underset{x \in \mathbb{R}^n}{\text{minimize}} \quad \max\left\{\underbrace{\tfrac{1}{2}x^T P x + p^T x + \alpha}_{\text{denoted } \psi_1(x)}, \, 0\right\} - \left(\underbrace{\tfrac{1}{2}x^T Q x + q^T x + \beta}_{\text{denoted } \psi_2(x)}\right) \tag{5}$$

has finitely many d-stationary values, where $P$ and $Q$ are symmetric matrices in $\mathbb{R}^{n \times n}$, $p$ and $q$ are $n$-vectors, and $\alpha$ and $\beta$ are scalars. By Corollary 5(a), we may assume that $P \neq 0$.

Let $\bar{x}$ be a d-stationary point of (5). If $\psi_1(\bar{x}) \neq 0$, then locally $\psi_1(x) \neq 0$ near $\bar{x}$. Thus $\bar{x}$ is a F-differentiable point of $\psi_1$. By Corollary 5(b), it therefore suffices to consider the case where $\psi_1(\bar{x}) = 0$. By the d-stationarity of $\bar{x}$, we have

$$\max\left((P\bar{x} + p)^T d, \, 0\right) - (Q\bar{x} + q)^T d \geq 0, \quad \forall d \in \mathbb{R}^n.$$

The latter condition is equivalent to the existence of a scalar $\bar{\lambda} \in [0, 1]$ such that $\bar{\lambda}(P\bar{x}+p) - (Q\bar{x}+q) = 0$, i.e.,

$$(\bar{\lambda} P - Q)\bar{x} = -(\bar{\lambda} p - q). \tag{6}$$

The argument below aims at showing that there are only finitely many values of $\psi_2(\bar{x})$ with $\bar{x}$ satisfying (6) and $\psi_1(\bar{x}) = 0$. We first notice that there are only finitely many real scalars $\lambda$ such that $\lambda P - Q$ is not invertible, since $\det(\lambda P - Q) = 0$ is a polynomial of $\lambda$ of degree $n$ (notice that $P \neq 0$). The remaining proof is divided into several cases.

• **Case 1:** $\bar{\lambda} \in [0, 1]$ is one of finitely many values for which $\det(\lambda P - Q) = 0$. The equation (6) shows that $\bar{x}$ is a stationary point of the (differentiable) quadratic function:

$$q_{\bar{\lambda}}(x) \triangleq \tfrac{1}{2}x^T\left[\bar{\lambda} P - Q\right]x + x^T\left(\bar{\lambda}p - q\right) + \bar{\lambda}\alpha - \beta = \bar{\lambda}\psi_1(x) - \psi_2(x).$$

Since for each fixed $\lambda$, being a quadratic function, $q_\lambda$ has at most one stationary value, by noticing that $q_{\bar{\lambda}}(\bar{x}) = -\psi_2(\bar{x})$, it follows that $\psi_2(\bar{x})$ is one of finitely many values in this case.

• **Case 2:** $\bar{\lambda}$ is such that $\det(\bar{\lambda}P - Q) \neq 0$. Without loss of generality, we may write all the distinct real solutions within the interval $[0, 1]$ of this determinantal equation as $\lambda_1 < \lambda_2 < \ldots < \lambda_m$ for some positive integer $m$, if at least one exists. For notational convenience, we let $\lambda_0 \triangleq 0 \leq \lambda_1$ and $\lambda_{m+1} \triangleq 1 \geq \lambda_m$. The discussion below thus also includes the case that $\det(\lambda P - Q) = 0$ has no solution in $(0, 1)$. For any $i = 0, 1, \cdots, m$ and any $\lambda \in (\lambda_i, \lambda_{i+1})$, the matrix $\lambda P - Q$ is invertible. In what follows, we may assume without loss of generality that the scalar $\bar{\lambda}$ satisfying (6) is such that $\bar{\lambda} \in (\lambda_i, \lambda_{i+1})$ for some $i \in \{0, 1, \cdots, m\}$. Thus, we obtain $\bar{x} = -(\bar{\lambda}P - Q)^{-1}(\bar{\lambda}p - q)$. Consider the vector function

$x(\lambda) \triangleq (\lambda P - Q)^{-1}(\lambda p - q)$ for $\lambda$ in $(\lambda_i, \lambda_{i+1})$. Either this function is identically equal to a constant, in which case $\bar{x}$ is a constant vector independent of $\lambda$, implying that $\psi(\bar{x})$ is one of finitely many values because there are only $(m+1)$ such intervals. Otherwise, by Lemma 6, $x(\lambda)$ is a nonzero rational vector-valued function of $\lambda$ that can be written as

$$x(\lambda) = \frac{\sum_{k=1}^{n} \lambda^k c^k + c^0}{(\lambda - \lambda_1)(\lambda - \lambda_2) \cdots (\lambda - \lambda_m) \left( \sum_{k=1}^{n-m} \gamma_k \lambda^k + \gamma_0 \right)}, \quad \lambda \in (\lambda_i, \lambda_{i+1})$$

for some vectors $\{c^k\}_{k=0}^{n} \in \mathbb{R}^n$ and some scalars $\{\gamma_k\}_{k=0}^{n-m}$ satisfying $\sum_{k=1}^{n-m} \gamma_k \lambda^k + \gamma_0 \neq 0$ for any $\lambda \in (\lambda_i, \lambda_{i+1})$. Consider the rational function $\psi_1(x(\lambda))$ for $\lambda \in (\lambda_i, \lambda_{i+1})$ and note that $\psi_1(x(\bar{\lambda})) = 0$. If this function is not identically equal to zero in this interval, then being a rational function in $\lambda$, it has only finitely many zeros in the interval, thus there are only finitely many $\bar{\lambda}$. We can now employ the argument in the previous case. Otherwise, if $\psi_1(x(\lambda)) \equiv 0$ for all $\lambda \in (\lambda_i, \lambda_{i+1})$, then

$$\frac{d\psi_1(x(\lambda))}{d\lambda} = (Px(\lambda) + p)^T \frac{dx(\lambda)}{d\lambda} = 0, \quad \forall \lambda \in (\lambda_i, \lambda_{i+1}),$$

which yields

$$\frac{d\psi_2(x(\lambda))}{d\lambda} = (Qx(\lambda) + q)^T \frac{dx(\lambda)}{d\lambda} = \lambda(Px(\lambda) + p)^T \frac{dx(\lambda)}{d\lambda} = 0, \quad \forall \lambda \in (\lambda_i, \lambda_{i+1}),$$

implying that $\psi_2(x(\lambda))$ is a constant in the interval $(\lambda_i, \lambda_{i+1})$. Thus, in either case, there are only finitely many possible values for $\psi_2(\bar{x})$ in the interval $(\lambda_i, \lambda_{i+1})$.

The proof can now be completed by summarizing the different cases considered above, in each of which $\psi_2(\bar{x})$ is one of finitely many values as claimed. $\square$

Although the above result does not impose any restriction on the number of active pieces of the second "max" of $\psi$, it pertains to d-stationary points with at most two active pieces in the first "max" function. At this time, we are not able to extend the arguments to three or more such active pieces without some restrictions. The difficulty with extending this proof lies in the proper treatment of the zeros of multivariate polynomial functions instead of a univariate polynomial. The next result allows an arbitrary number of active pieces in two special cases. Notice that the function $\psi$ in the proposition is a dc function that is neither convex nor concave. Part (a) of the proposition is a special case of Proposition 2(b); we include it in order to contrast with part (b) which requires a linear independence assumption of gradients.

**Proposition 9.** Suppose that $X$ is a polyhedral set. Let $\psi(x)$ be given by (2) where the functions $\psi_{ji}$ are all quadratic and satisfy one of the following two conditions:

(a) each $\psi_{1i}$ is convex for $i = 1, \cdots k_1$ and each $\psi_{2i}$ is concave for $i = 1, \cdots, k_2$; or

(b) each $\psi_{i1}$ is concave for $i = 1, \cdots k_1$, each $\psi_{2i}$ is strictly convex for $i = 1, \cdots, k_2$, and the gradients $\{\nabla \psi_{1i}(\bar{x})\}_{i \in \mathcal{A}_1(\bar{x})}$ are linearly independent for all $\bar{x} \in \mathcal{D}_{(\psi, X)}$.

Then the set $\psi\left(\mathcal{D}_{(\psi, X)}\right)$ is finite under (a) and the set $\mathcal{D}_{(\psi, X)}$ is finite under (b).

**Proof.** It suffices to prove statement (b). Without loss of generality, we may assume that $X = \mathbb{R}^n$. Moreover, by Lemma 1, it suffices to prove the result for

$$\psi(x) = \max_{1 \leq i \leq k_1} \psi_{1i}(x) - q_2(x),$$

where $q_2$ is a strictly convex quadratic function. For $i = 1, \cdots, k_1$, we may write $\psi_{1i}(x) \triangleq \frac{1}{2}x^T Q^i x + x^T q^i + \alpha_i$ for some symmetric negative semidefinite matrix $Q^i$, vector $q^i$, and scalar $\alpha_i$; we also write $q_2(x) \triangleq \frac{1}{2}x^T P x + p^T x + \beta$ for some symmetric positive definite matrix $P$, vector $p$, and scalar $\beta$. Let $\bar{x} \in \mathcal{D}_{(\psi,X)}$ be arbitrary. We note that the condition of linear independent gradients of the argmax functions remain valid for any of its sub-family of functions, at the same d-stationary points. As such, by induction on $k_1$, we may assume without loss of generality that $\psi_{1i}(\bar{x}) = 0$ for all $i = 1, \cdots, k_1$; i.e., $\mathcal{A}_1(\bar{x}) = \{1, \cdots, k_1\}$. Since $X = \mathbb{R}^n$, there exists a family of nonnegative scalars $\{\bar{\lambda}_i\}_{i=1}^{k_1}$ satisfying $\sum_{i=1}^{k_1} \bar{\lambda}_i = 1$ such that

$$P\bar{x} + p = \sum_{i=1}^{k_1} \bar{\lambda}_i \left( Q^i \bar{x} + q^i \right). \tag{7}$$

Since each $Q^i$ is negative semidefinite, the matrix $\bar{Q} \triangleq P - \sum_{i=1}^{k_1} \bar{\lambda}_i Q^i$ is positive definite. For an element $\lambda$ of the set $\Lambda \triangleq \left\{ \lambda \in \mathbb{R}_+^{k_1} \mid \sum_{i=1}^{k_1} \lambda_i = 1 \right\}$, let $Q(\lambda) \triangleq P - \sum_{i=1}^{k_1} \lambda_i Q^i$ which is positive definite; also let $x(\lambda) = Q(\lambda)^{-1} \left( \sum_{i=1}^{k_1} \lambda_i q^i - p \right)$. We claim that every zero $\widehat{\lambda} \in \Lambda$ of the composite mapping $F : \Lambda \to \mathbb{R}^{k_1}$ defined by

$$F(\lambda) \triangleq \begin{pmatrix} (\psi_{11} - \psi_{1k_1}) \circ x(\lambda) \\ \vdots \\ (\psi_{1k_1-1} - \psi_{1k_1}) \circ x(\lambda) \\ \sum_{i=1}^{k_1} \lambda_i - 1 \end{pmatrix}, \quad \text{for } \lambda \in \Lambda,$$

is isolated. For this, it suffices to show that the Jacobian matrix $JF(\widehat{\lambda})$ is nonsingular. Since by definition, $Q(\lambda) x(\lambda) = \sum_{i=1}^{k_1} \lambda_i q^i - p$, it follows that

$$\frac{\partial x(\lambda)}{\partial \lambda_j} = Q(\lambda)^{-1} \nabla \psi_{1j}(x(\lambda)), \quad j = 1, \cdots, k_1.$$

Hence by the chain rule, we deduce that

$$\frac{\partial [\psi_{1i} - \psi_{1k_1}] \circ x(\lambda)}{\partial \lambda_j} = [\nabla \psi_{1i}(x(\lambda)) - \nabla \psi_{1k_1}(x(\lambda))]^T \frac{\partial x(\lambda)}{\partial \lambda_j}$$
$$= [\nabla \psi_{1i}(x(\lambda)) - \nabla \psi_{1k_1}(x(\lambda))]^T Q(\lambda)^{-1} \nabla \psi_{1j}(x(\lambda)), \quad \forall i = 1, \cdots, k_1 - 1 \text{ and } j = 1, \cdots k_1.$$

Consequently,

$$JF(\lambda) = \begin{bmatrix} \widehat{\Psi}(\lambda)^T Q(\lambda)^{-1} \Psi(\lambda) \\ \mathbf{1}_{k_1}^T \end{bmatrix} \in \mathbb{R}^{k_1 \times k_1},$$

where $\widehat{\Psi}(\lambda)$ is the $n \times (k_1 - 1)$ matrix with the columns being $\nabla \psi_{1i}(x(\lambda)) - \nabla \psi_{1k_1}(x(\lambda))$ for $i = 1, \cdots, k_1 - 1$, $\Psi(\lambda)$ is the $n \times k_1$ matrix with the columns being $\nabla \psi_{1i}(x(\lambda))$ for $i = 1, \cdots, k_1$, and $\mathbf{1}_{k_1}$ is the $k_1$-dimensional vector of all ones. If $\widehat{\lambda} \in \Lambda$ is such that $F(\widehat{\lambda}) = 0$, then $x(\widehat{\lambda})$ belongs to $\mathcal{D}_{(\psi,X)}$. Hence the matrix $\Psi(\widehat{\lambda})$ has linearly independent columns by assumption; thus so does the matrix $\widehat{\Psi}(\widehat{\lambda})$. To show that $JF(\widehat{\lambda})$ is nonsingular, let $w \in \mathbb{R}^{k_1}$ be such that $JF(\widehat{\lambda})w = 0$. We then have $w_{k_1} = -\sum_{i=1}^{k_1-1} w_i$. Hence we deduce
$$0 = \widehat{\Psi}(\widehat{\lambda})^T Q(\widehat{\lambda})^{-1} \Psi(\widehat{\lambda})w = \widehat{\Psi}(\widehat{\lambda})^T Q(\widehat{\lambda})^{-1} \widehat{\Psi}(\widehat{\lambda})\widehat{w},$$
where $\widehat{w}$ is the first $(k_1 - 1)$ components of $w$. Since the matrix $\widehat{\Psi}(\widehat{\lambda})$ has linearly independent columns and $Q(\widehat{\lambda})^{-1}$ is positive definite, it follows that $\widehat{w} = 0$; hence $w_{k_1} = 0$. Consequently, the Jacobian matrix $JF(\widehat{\lambda})$ is nonsingular, establishing the isolatedness of $\widehat{\lambda}$. Having only isolated zeros in the compact set $\Lambda$ with $\bar{\lambda}$ being one of them, the function $F$ has only finitely many zeros in $\Lambda$. It therefore follows that $\bar{\lambda}$ is one of finitely many values; hence so is $\bar{x} = x(\bar{\lambda})$. □

Under the assumptions in part (b) of the above proposition, the set $\mathcal{D}_{(\psi,X)}$ of d-stationary points is finite. Thus the condition of linear independence of gradients seems stronger than necessary for the finiteness of the d-stationary values and thus deserves some discussion. This condition can be equivalently phrased in a constraint formulation of the first pointwise-max term in the objective function $\psi$. Indeed setting aside the second pointwise-max term, we may cast the problem (3) as $k_1$ problems in which the first pointwise-max function is transformed into constraints:

$$\left\{ \begin{array}{ll} \underset{x \in X}{\text{minimize}} & \psi_{1i}(x) - \max_{1 \leq i \leq k_2} \psi_{2i}(x) \\ \text{subject to} & \psi_{1i}(x) \geq \psi_{1j}(x), \quad \text{for all } j \neq i \end{array} \right\}_{i=1}^{k_1}.$$

Due to the nonlinearity of the constraints, some qualifications are needed to understand the tangent vectors of the feasible set, which are key to the definition of stationarity conditions of a non-convex constrained optimization problem in general, and each of the above $k_1$ quadratically constrained problems in particular; see Section 4. Among these constraint qualifications (CQs), the linear independence of the gradients of the active constraints (LICQ) is the most prominent one. We have the following linkage between this CQ and the LI condition in part (b) of Proposition 9. The proof is a trivial exercise about linear independence and thus omitted.

**Proposition 10.** Let $\{g_i\}_{i=1}^k$ be a family of continuously differentiable functions defined on an open set $\Omega \subseteq \mathbb{R}^n$. The gradients $\{\nabla g_i(\bar{x})\}_{i \in \mathcal{A}(\bar{x})}$ where $\mathcal{A}(\bar{x}) \triangleq \underset{1 \leq i \leq k}{\text{argmax}}\, g_i(\bar{x})$, are linearly independent if and only if for every $i \in \{1, \cdots, k\}$, the LICQ holds for the constraints $\{g_i(x) \geq g_j(x)\}_{j \neq i}$ at $\bar{x}$.

Regrettably, while some of them are not easy to prove, the collection of the obtained results has not fully resolved the following main question of this paper, which we pose as a conjecture for future research:

**Conjecture:** A piecewise quadratic function of the difference-max type (2) has only finitely many d-stationary values on polyhedra.

### 3.1 Composite piecewise extensions.

Based on the established results, we may further show the finite number of d-stationary values for a class of composite piecewise programs. Consider the composite program

$$\underset{x \in X}{\text{minimize}} \quad f(x) \triangleq \varphi \circ \psi(x) \tag{8}$$

where $\varphi : \mathbb{R} \to \mathbb{R}$ is a univariate function and $X \subseteq \mathbb{R}^n$ is a convex polyhedral set. The following proposition generalizes some of the propositions in the previous subsection.

**Proposition 11.** The problem (8) has only finitely many d-stationary values if $\varphi$ is a convex or concave univariate function, and $\psi$ satisfies one of the following four assumptions:

(a) $\psi$ is a piecewise linear-quadratic function;

(b) $\psi(x) = \psi_1(x) + \max\limits_{1 \leq i \leq k_2} \psi_{2i}(x)$ with $\psi_1$ being a piecewise affine function and each $\psi_{2i}$ being a convex quadratic function;

(c) $\psi(x) = \max(\psi_{11}(x), \psi_{12}(x)) - \max(\psi_{21}(x), \psi_{22}(x))$ with each $\psi_{ij}$ being a quadratic function.

(d) $\psi$ is given by (2) with the two families of quadratic functions $\{\psi_{1i}\}_{i=1}^{k_1}$ and $\{\psi_{2i}\}_{i=1}^{k_2}$ satisfying the conditions in part (b) of Proposition 9.

**Proof.** Let $\bar{x} \in X$ be a d-stationary point of $f$. Then by the chain rule of directional differentiation,

$$0 \leq f'(\bar{x}; x - \bar{x}) = \varphi'(\psi(\bar{x}); \psi'(\bar{x}; x - \bar{x})) = |\psi'(\bar{x}; x - \bar{x})| \varphi'(\psi(\bar{x}); \pm 1), \quad \forall\, x \in X,$$

where the $\pm 1$ depends on the sign of the directional derivative $\psi'(\bar{x}; x - \bar{x})$. [If this derivative is zero, then $\pm 1$ is immaterial.] If $\varphi$ is a convex univariate function, one of the following three cases must hold: (i) $\varphi'(\psi(\bar{x}); \pm 1) \geq 0$; (ii) $\varphi'(\psi(\bar{x}); +1) \geq 0 > \varphi'(\psi(\bar{x}); -1)$; or (iii) $\varphi'(\psi(\bar{x}); +1) < 0 \leq \varphi'(\psi(\bar{x}); -1)$. If case (i) holds, then $\psi(\bar{x}) \in \underset{t \in \mathbb{R}}{\operatorname{argmin}}\, \varphi(t)$ and thus, $f(\bar{x}) = \underset{x \in X}{\operatorname{argmin}}\, f(x)$, i.e., $\bar{x}$ is the global minimizer of $f$. If case (ii) holds, we have $\psi'(\bar{x}; x - \bar{x}) \geq 0$ for all $x \in X$. Then $\bar{x}$ is a d-stationary point of $\psi$ on $X$. Similarly for case (iii), we have $(-\psi)'(\bar{x}; x - \bar{x}) \geq 0$ for all $x \in X$, i.e., $\bar{x}$ is a d-stationary point of $-\psi$ on $X$. It follows from Propositions 2, 8, and 9 that if $\psi$ satisfies one of assumptions (a)–(d), then both $\psi(\mathcal{D}_{(\psi,X)})$ and $(-\psi)(\mathcal{D}_{(-\psi,X)})$, and thus $\psi(\mathcal{D}_{(-\psi,X)})$ too, are finite sets; hence $f$ has a finite number of d-stationary values in cases (ii) and (iii).

If $\varphi$ is a concave function, one of the following three cases must hold: (i) $\varphi'(\psi(\bar{x}); \pm 1) \leq 0$; (ii) $\varphi'(\psi(\bar{x}); +1) > 0 \geq \varphi'(\psi(\bar{x}); -1)$; or (iii) $\varphi'(\psi(\bar{x}); +1) \leq 0 < \varphi'(\psi(\bar{x}); -1)$. By similar arguments we can establish the desired results. $\square$

**Remark 12.** Although Proposition 11 holds for $\varphi$ being either a convex function or a concave function, it cannot be generalized to a dc function $\varphi$. In fact, a univariate dc function, even if it is continuously differentiable, may have infinitely many d-stationary values. This can be seen from the example $\varphi(t) = \varphi_1(t) - \varphi_2(t)$, where $t \in \mathbb{R}$,

$$\varphi_1(t) = \begin{cases} \dfrac{3}{2}t^2 - 2nt + 2n^2 + n + 1 & \text{if } t \in [2n, 2n+1) \\ \dfrac{1}{2}t^2 + 2(n+1)t - 2n^2 - 3n & \text{if } t \in [2n+1, 2n+2), \end{cases} \quad n = 0, \pm 1, \pm 2, \cdots$$

and

$$\varphi_2(t) = \begin{cases} t^2 + 2nt - 2n^2 - n - 1 & \text{if } t \in [2n, 2n+1) \\ \dfrac{3}{2}t^2 - 2(n+1)t + 2n^2 + 3n & \text{if } t \in [2n+1, 2n+2), \end{cases} \quad n = 0, \pm 1, \pm 2, \cdots.$$

One may verify that

$$\varphi_1'(t) = \begin{cases} 3t - 2n & \text{if } t \in [2n, 2n+1) \\ t + 2(n+1) & \text{if } t \in [2n+1, 2n+2), \end{cases} \qquad \varphi_2'(t) = \begin{cases} t + 2n & \text{if } t \in [2n, 2n+1) \\ 3t - 2(n+1) & \text{if } t \in [2n+1, 2n+2). \end{cases}$$

Since both derivatives are increasing functions on $\mathbb{R}$, it follows that $\varphi_1$ and $\varphi_2$ are convex functions. Moreover, we have

$$\varphi'(t) = \begin{cases} 2t - 4n & \text{if } t \in [2n, 2n+1) \\ -2t + 4(n+1) & \text{if } t \in [2n+1, 2n+2), \end{cases} \quad n = 0, \pm 1, \pm 2, \cdots,$$

i.e., $\varphi'(t) \geq 0$ and $\varphi'(t) = 0$ if and only if $t = 2n$, which further implies that $\varphi$ is a strictly increasing function on $\mathbb{R}$. Therefore, $\varphi$ has infinite many d-stationary points $t = 2n$ and infinite many d-stationary values $2n+2$ for any integer $n$; see Figure 2 for an illustration of the function $\varphi$ and its derivative. Notice that the function $\varphi$ is not piecewise linear-quadratic since there are infinitely many different quadratic pieces. $\square$

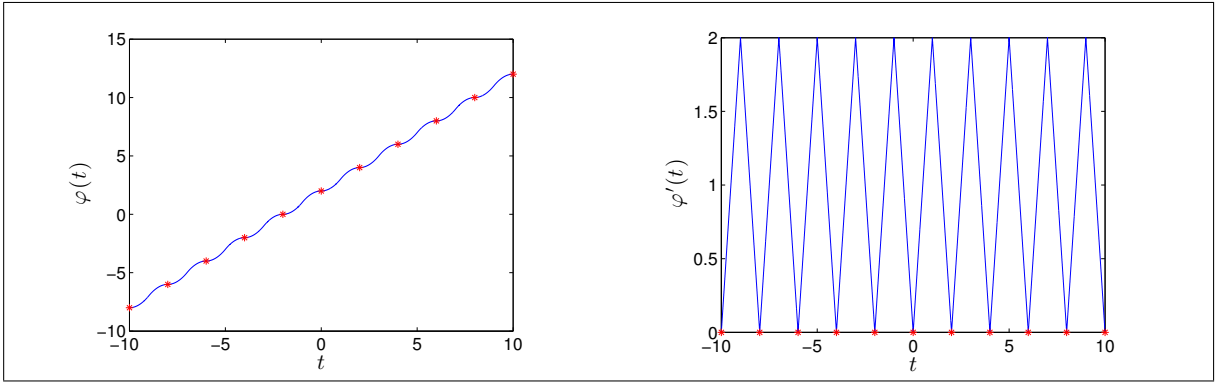

Figure 2: the dc function $\varphi$ and its derivative

## 4 A Doubly Quadratic Inequality Constrained Piecewise Program.

Consider a class of doubly quadratic inequality constrained piecewise quadratic programs as follows:

$$\begin{aligned}
\underset{x \in \mathbb{R}^n}{\text{minimize}} \quad & \psi(x) \triangleq \underbrace{\tfrac{1}{2} x^T P^0 x + (p^0)^T x}_{\text{denoted } q_1(x)} - \max_{1 \leq i \leq k_2} \left\{ \underbrace{\tfrac{1}{2} x^T P^i x + (p^i)^T x + \alpha_i}_{\text{denoted } q_{2i}(x)} \right\} \\
\text{subject to} \quad & \beta_1 \leq \tfrac{1}{2} x^T Q x + c^T x \leq \beta_2, \quad \text{and} \quad Ax \leq b,
\end{aligned} \quad (9)$$

where $P^0$, $P^i$, and $Q \in \mathbb{R}^{n \times n}$ are symmetric matrices, $A$ is an $m \times n$ matrix, $p^0$, $p^i$, and $c$ are $n$-vectors, $b$ is an $m$-vector, $\alpha_i$, $\beta_1$, and $\beta_2$ are scalars with $\beta_1 \leq \beta_2$. The discussion below allows for the possibilities that either $\beta_1 = -\infty$ or $\beta_2 = \infty$, which corresponds to the case that effectively there is only one quadratic inequality constraint in the problem. The case $\beta_1 = \beta_2$ is also included; this corresponds to the case where the two quadratic inequalities collapse into one quadratic equality. The program (9) arises frequently in various optimization problems. Without the pointwise-max term in the objective, it arises, for example, in trust region methods for solving (equality constrained) nonlinear programs [13, 25] with positive definite scaling matrix $Q$:

$$\begin{aligned}
\underset{x \in \mathbb{R}^n}{\text{minimize}} \quad & \psi(x) \triangleq \tfrac{1}{2} x^T P^0 x + (p^0)^T x \\
\text{subject to} \quad & \beta_1 \leq \tfrac{1}{2} x^T Q x \leq \beta_2, \quad Ax = b,
\end{aligned} \quad (10)$$

and the perturbed eigenvalue problems [15, 38]:

$$\begin{aligned}\underset{x\in\mathbb{R}^n}{\text{minimize}} \quad & x^T P^0 x - (p^0)^T x \\ \text{subject to} \quad & \|x\|_2 = \beta \quad (\Leftrightarrow \beta^2 \leq x^T x \leq \beta^2)\end{aligned}$$

Let $X$ denote the feasible set of the problem (9). Since $X$ is not convex in general, we need to extend the d-stationarity, which was defined for a convex feasible set to that of a B(ouligand)-stationary point [28]. The Bouligand tangent cone of $X$ at $x \in X$, denoted $\mathcal{T}(x; X)$, is a closed cone (not necessarily convex) whose elements are vectors $d \in \mathbb{R}^n$ for which a sequence $x^k \subseteq X$ of vectors converging to $x$ and a sequence of positive scalars $\{\tau_k\}$ converging to zero exist such that

$$d = \lim_{k\to\infty} \frac{x^k - x}{\tau_k}.$$

A vector $\bar{x} \in X$ is a B-stationary point of (9) if it satisfies $\psi'(\bar{x}; d) \geq 0$ for all $d \in \mathcal{T}(\bar{x}; X)$; or equivalently, for any $\bar{i} \in \mathcal{A}_2(\bar{x}) \triangleq \underset{1\leq i \leq k_2}{\operatorname{argmax}}\ q_{2i}(\bar{x})$, it holds

$$[\nabla q_1(\bar{x}) - \nabla q_{2\bar{i}}(\bar{x})]^T d \geq 0, \quad \forall\, d \in \mathcal{T}(\bar{x}; X). \tag{11}$$

When the set $X$ is convex, B-stationarity reduces to d-stationarity. The above condition is a primal description of a B-stationary point. Constraint qualifications (CQs) are needed for an analytical description of the tangent cone, which leads to a primal-dual description of the stationary point that involves the constraint multipliers. One of the most general CQs is Abadie's CQ, implying the equality between the tangent cone $\mathcal{T}(\bar{x}; X)$ and the linearization cone of $X$ at $\bar{x}$, i.e.,

$$\mathcal{T}(\bar{x}; X) = \left\{ v \in \mathbb{R}^n \mid \begin{array}{ll} v^T(Q\bar{x} + c) \geq 0 & \text{if } \tfrac{1}{2}\bar{x}^T Q\bar{x} + c^T \bar{x} = \beta_1 \\ v^T(Q\bar{x} + c) \leq 0 & \text{if } \tfrac{1}{2}\bar{x}^T Q\bar{x} + c^T \bar{x} = \beta_2 \\ A_{i\bullet} v \leq 0 & \text{for all } i \in \mathcal{I}(\bar{x}) \triangleq \{ i \mid A_{i\bullet}\bar{x} = b_i \} \end{array} \right\}.$$

Under this condition, (11) holds for a given $\bar{i} \in \mathcal{A}_2(\bar{x})$ if and only if there exist multipliers $\lambda_j^{\bar{i}}$ for $j = 1, 2$, and $\mu^{\bar{i}} \in \mathbb{R}^m$ such that

$$\nabla q_1(\bar{x}) - \nabla q_{2\bar{i}}(\bar{x}) + (\lambda_2^{\bar{i}} - \lambda_1^{\bar{i}})(Q\bar{x} + c) + A^T \mu^{\bar{i}} = 0$$

$$0 \leq \lambda_1^{\bar{i}} \perp \tfrac{1}{2}\bar{x}^T Q\bar{x} + c^T\bar{x} - \beta_1 \geq 0$$

$$0 \leq \lambda_2^{\bar{i}} \perp \beta_2 - \left(\tfrac{1}{2}\bar{x}^T Q\bar{x} + c^T\bar{x}\right) \geq 0$$

$$0 \leq \mu^{\bar{i}} \perp b - A\bar{x} \geq 0,$$

where the $\perp$ notation expresses the complementarity between the constraint slacks and associated multipliers. For any B-stationary point $\bar{x} \in X$ of (9), we call $\psi(\bar{x})$ a B-stationary value. With similar arguments to the proof of Proposition 8, we can show that the problem (9) has only finitely many B-stationary values.

**Proposition 13.** The problem (9) has finitely many B-stationary values under either one of the following two assumptions:

(a) Abadie's CQ holds at all B-stationary points;

(b) One of the two quadratic inequalities cannot be satisfied strictly on the polyehdron: $Z \triangleq \{x \in \mathbb{R}^n \mid Ax \leq b\}$.

**Proof.** We first prove the proposition under assumption (a). It suffices to show that $\psi$ has finitely many different values over the points $x$ satisfying (11) for a fixed $\bar{i} \in \{1,\ldots,k_2\}$. Similarly to the proof of Proposition 8, we may take $X$ to be the entire space $\mathbb{R}^n$. To proceed, we write $\widehat{P}^{\bar{i}} \triangleq P^0 - P^{\bar{i}}$ and $\widehat{p}^{\bar{i}} \triangleq p^0 - p^{\bar{i}}$. If $\beta_1 < \frac{1}{2}\bar{x}^T Q \bar{x} + c^T \bar{x} < \beta_2$, then the corresponding multipliers $\lambda_1^{\bar{i}} = \lambda_2^{\bar{i}} = 0$ and $\bar{x}$ is a stationary point of the quadratic program: $\displaystyle\minimize_{x \in \mathbb{R}^n} \frac{1}{2} x^T \widehat{P}^{\bar{i}} x + (\widehat{p}^{\bar{i}})^T x$, which has only finitely many stationary values. Otherwise if $\beta_1 < \frac{1}{2}\bar{x}^T Q\bar{x} + c^T \bar{x} = \beta_2$, then there exists $\bar{\lambda} \geq 0$ such that

$$\widehat{P}^{\bar{i}}\bar{x} + \widehat{p}^{\bar{i}} + \bar{\lambda}(Q\bar{x} + c) = 0. \tag{12}$$

Notice that the proof in Proposition 8 for the finite number of d-stationary values of (5) in fact does not rely on the condition that $\bar{\lambda} \in [0,1]$ in (6). Thus, by similar arguments to Proposition 8, we may show that $\frac{1}{2} x^T \widehat{P}^{\bar{i}} x + (\widehat{p}^{\bar{i}})^T x$ has finitely many values on $\bar{x}$ satisfying both (12) and the equality $\frac{1}{2}\bar{x}^T Q \bar{x} + c^T \bar{x} - \beta_2 = 0$. For the case that $\beta_1 = \frac{1}{2}\bar{x}^T Q\bar{x} + c^T \bar{x} = \beta_2$, there exists $\lambda \in \mathbb{R}$ such that (12) holds. We can thus also show that $\frac{1}{2} x^T \widehat{P}^{\bar{i}} x + (\widehat{p}^{\bar{i}})^T x$ has finitely many values on such $\bar{x}$.

Under assumption (b), suppose that the inequality $\frac{1}{2} x^T Q x + c^T x \leq \beta_2$ cannot be satisfied strictly on the polyhedron $Z$. In this case, the feasible set of (9), if nonempty, is equal to the set

$$\{\, x \in Z \mid \tfrac{1}{2} x^T Q x + c^T x = \beta_2 \,\} \;=\; \argmin\{\, \tfrac{1}{2} x^T Q x + c^T x \mid x \in Z \,\}.$$

In turn, if $x$ is an element of the latter argmin set, then there exists a multipier $\mu$ such that $(x,\mu)$ satisfies the mixed complementarity conditions:

$$Qx + c + A^T \mu = 0$$
$$0 \leq \mu \perp b - Ax \geq 0.$$

The set of pairs $(x,\mu)$ satisfying the latter conditions is the union of finitely many polyhedra, each of the form:
$$\begin{aligned} Qx + c + A^T\mu &= 0 \\ \mu_I &\geq 0 = (b - Ax)_I \\ \mu_J &= 0 \leq (b - Ax)_J, \end{aligned} \tag{13}$$

for some complementary pair of index subsets $I$ and $J$ of $\{1,\cdots,m\}$. Let $P_{IJ}$ be the polyehdron consisting of vectors $x$ for which there exists $\mu$ such that the pair $(x,\mu)$ satisfies (13). The set of d-stationary values of the quadratic function $q_1 - q_{2\bar{i}}$ is finite on each nonempty $P_{IJ}$. The proof will be completed if we can show that any d-stationary point of $q_1 - q_{2\bar{i}}$ on the feasible set of (9) is a d-stationary point of the same function on some such set $P_{IJ}$. Let $\bar{x}$ be a d-stationary point of $q_1 - q_{2\bar{i}}$ on the feasible set of (9). Let $(\bar{I}, \bar{J})$ be a pair of index sets such $\bar{x} \in P_{\bar{I}\bar{J}}$. We claim that $\bar{x}$ is a d-stationary point of $q_1 - q_{2\bar{i}}$ on $P_{\bar{I}\bar{J}}$. It suffices to show that $x - \bar{x}$ is a tangent vector of the feasible set of (9) at $\bar{x}$ for any $x \in P_{\bar{I}\bar{J}}$. We show this by verifying that $\frac{1}{2}(x^\tau)^T Q x^\tau + c^T x^\tau = \frac{1}{2}\bar{x}^T Q \bar{x} + c^T \bar{x}$ for all $\tau \in [0,1]$, where $x^\tau \triangleq \bar{x} + \tau(x - \bar{x})$. Let $\bar{\mu}$ and $\mu$ be such that the pairs $(\bar{x}, \bar{\mu})$ and $(x, \mu)$ both satisfy (13) corresponding to the pair $(\bar{I}, \bar{J})$. We have

$$\begin{aligned} \tfrac{1}{2}(x^\tau)^T Q x^\tau + c^T x^\tau &= \tfrac{1}{2}\bar{x}^T Q \bar{x} + c^T \bar{x} + \tau(Q\bar{x} + c)^T(x - \bar{x}) + \tfrac{\tau^2}{2}(x-\bar{x})^T Q(x - \bar{x}) \\ &= \tfrac{1}{2}\bar{x}^T Q \bar{x} + c^T \bar{x} - \tau(x - \bar{x})^T A^T \bar{\mu} - \tfrac{\tau^2}{2}(x-\bar{x})^T A^T(\mu - \bar{\mu}) \\ &= \tfrac{1}{2}\bar{x}^T Q \bar{x} + c^T \bar{x} \end{aligned}$$

because $(x - \bar{x})^T A^T \bar{\mu} = [(Ax - A\bar{x})_{\bar{I}}]^T \bar{\mu}_{\bar{I}} = 0 = (x - \bar{x})^T A^T (\mu - \bar{\mu})$. $\square$

It would be worthwhile to relate Proposition 13 to some existing results in the literature and comment on how the former result is different from the latter. First, for the trust region problem (10), if it is feasible and a so-called simultaneous diagonalization condition holds, i.e., there exists a scalar $\gamma$ such that $P^0 + \gamma Q$ is positive definite, then the problem is equivalent to a convex minimization problem with simple linear constraints [5]. The simultaneous diagonalization condition holds if $Q = I$, as in the case of the perturbed eigenvalues problems. In addition, a globally optimal solution in special instances of (9) may be found in polynomial time via semidefinite programming (SDP) relaxation if the objective is given by a quadratic function (i.e., without the pointwise max term); these instances are: (i) there is a single quadratic inequality constraint and no linear (equality or inequality) constraints; (ii) there is a single convex quadratic inequality constraint and a single linear inequality constraint; (iii) the quadratic functions in both the objective and the constraints are homogeneous quadratic functions (i.e., no linear term in the quadratic function) and there is no linear constraint (in fact, two quadratic inequality constraints with different Hessian matrices are allowed for this case). The first two cases are studied in [39] while the last one is discussed in [31, 42]. See also [43, 19, 4, 18] for other conditions of the exact SDP relaxation for the quadratically constrained quadratic programs when no linear constraint appears in (9). Lastly, a recent paper [1] has studied a conic quadratically constrained quadratic program and shown that it can be "lifted" to an equivalent convex "completely positive program". Such an equivalence is in terms of the globally optimal solutions of the problems; there is no mention about stationary solutions of the original non-convex quadratically constrained problem.

For those cases where a globally optimal solution of (9) can be computed efficiently by convexification of some sort, Proposition 13 may not be interesting. However, admitting an arbitrary number of linear constraints and two quadratic inequalities (albeit with the same quadratic form), Proposition 13 encompasses a much broader class of piecewise quadratic programming problems than those studied in the above cited references. One clear distinction is the inclusion of the pointwise-max, thus nonsmooth piecewise quadratic term in the objective; another noteworthy point is that all the quadratic functions can be indefinite with no simultaneous diagonalization condition assumed. In contrast to the cited results all of which pertain to global optimality, Proposition 13 pertains to d-stationary values and makes no claim about the globally optimal solutions.

## 5 Applications of Results.

This section presents classes of optimization problems to which the results in the previous sections are applicable. At this time, these results should be considered as providing new theoretical properties of the problems as we make no claim about how the results can facilitate the numerical solutions of these non-convex problems, particularly for the computation of globally optimal solutions. Before discussing these problems, we give one consequence of the finiteness of B-stationary values pertaining to the sequential convergence of objective values. We state the following result for the general problem (3) without requiring the feasible set be convex.

**Proposition 14.** Let the function $\psi$ be continuous. Suppose that the set of B-stationary values of the optimization problem (3) is finite. If $\{x^k\}$ is a bounded sequence with the property that every one of its accumulation points is B-stationary for the problem, then the sequence of objective values $\{\psi(x^k)\}$ converges (to a B-stationary value) if and only if $\{\psi(x^{k+1}) - \psi(x^k)\} \to 0$.

**Proof.** It suffices to show the "if" statement. In turn, according to [10, Proposition 8.3.10], it suffices to show that the scalar sequence $\{\psi(x^k)\}$, which must be bounded, has only finitely many accumulation points. Let $\kappa$ be an infinite subset of $\{1, 2, \cdots\}$ such that the subsequence $\{\psi(x^k)\}_{k \in \kappa}$ converges to a value $\psi_\infty$. The subsequence of vectors $\{x^k\}_{k \in \kappa}$ must have an accumulation point, say $x^\infty$, which must be B-stationary. Thus $\psi_\infty = \psi(x^\infty)$ is a B-stationary value. Since there are only finitely many B-stationary values, the sequence $\{\psi(x^k)\}$ has (at least one and) only finitely many accumulation points, thus converges. □

The key point of the above result is that the sequence $\{\psi(x^k)\}$ of objective values is not required to be monotonic. This may occur in an iterative algorithm wherein the values of the original objective function may not be decreasing, but those of a merit function are decreasing. Such a situation happens to the alternating direction method of multipliers for solving a class of linearly constrained nonconvex problems [29]. We also illustrate this situation in Subsection 5.3 pertaining to a majorization-minimization method for solving the least-squares piecewise affine regression problem.

## 5.1 Mathematical programs with complementarity constraints.

Mathematical programs with complementarity constraints (MPCCs) are optimization problems subject to disjunctive constraints represented by the complementarity relations. It is a subclass of mathematical programs with equilibrium constraints; see [21, 11] for comprehensive discussions on this subject. We consider the convex programs with linear complementarity constraints:

$$\begin{aligned} \underset{x,y}{\text{minimize}} \quad & f(x,y) \\ \text{subject to} \quad & 0 \leq y \perp F(x,y) \triangleq q + Nx + My \geq 0, \quad \text{and} \quad (x,y) \in Z, \end{aligned} \quad (14)$$

where $f : \mathbb{R}^{n+m} \to \mathbb{R}$ is a convex function, $q \in \mathbb{R}^m$ is a vector, $N \in \mathbb{R}^{m \times n}$ and $M \in \mathbb{R}^{m \times m}$ are two matrices, and $Z$ is a polyhedral set.

It is known that the MPCC has no feasible point that satisfies the Mangasarian–Fromovitz constraint qualification. One way to solve the MPCC is via the penalty method that removes the complementarity constraint $y^T F(x,y) = 0$ from the above formulation by adding a penalty term $\gamma p(y, F(x,y))$ to the objective function for some positive penalty parameter $\gamma$. One particular choice of the penalty function is

$$p(y, F(x,y)) = \sum_{i=1}^{m} \min(y_i, F_i(x,y)) = \min_{I \subseteq \{1,\ldots,m\}} \left( \sum_{i \in I} y_i + \sum_{i \in I^c} F_i(x,y) \right),$$

where $I^c$ is the complement of $I$. Such a penalty function is the pointwise minimum of a finite number of affine functions, a special concave piecewise affine functions. The resulting penalized problem is

$$\begin{aligned} \underset{x,y}{\text{minimize}} \quad & f(x,y) - \gamma \max_{I \subseteq \{1,\ldots,m\}} \left( -\sum_{i \in I} y_i - \sum_{i \in I^c} F_i(x,y) \right) \\ \text{subject to} \quad & y \geq 0, \quad F(x,y) \geq 0, \quad (x,y) \in Z. \end{aligned}$$

Based on Proposition 2(b), the above nonlinear program has only finitely many d-stationary values.

A prominent subclass of MPCCs consists of quadratic programs with complementarity constraints, in which the objective function $f(x,y)$ in (14) is a quadratic function [2]. We may alternatively adopt the following quadratic penalty function for the complementarity constraint [17, 32]

$$p(x, F(x,y)) = x^T F(x,y),$$

which leads to a polyhedral constrained quadratic programming. By [22, Lemma 3.1], the penalized problem has finitely many d-stationary values. In fact, by Proposition 2, for the quadratic penalized problem to have finitely many d-stationary points, the objective function $f$ in (14) only needs to be a piecewise linear-quadratic function, or the pointwise max of finitely many quadratic functions.

Besides the penalty approach, another strategy to solve the MPCC is via the regularization method [35, 32], which approximates (14) by relaxing the complementarity constraint, resulting in the following nonlinear program parameterized by the positive tolerance $\varepsilon$:

$$\begin{aligned}&\underset{x,y}{\text{minimize}}\quad f(x,y)\\&\text{subject to}\quad y\geq 0,\quad q+Nx+My\geq 0,\quad (x,y)\in Z,\quad y^T(q+Nx+My)\leq\varepsilon.\end{aligned}$$

If the function $f$ is quadratic, the above problem is a singly quadratic inequality constrained quadratic programming. By Proposition 13, it has finitely many B-stationary values under Abadie's CQ.

## 5.2 Two-stage (indefinite) stochastic quadratic programs.

We consider a quadratic programming based recourse function [20, 27] in two-stage stochastic programming [6, 37], in which the first stage is also defined by a quadratic problem:

$$\begin{aligned}&\underset{x\in\mathbb{R}^n}{\text{minimize}}\quad \frac{1}{2}x^TPx+c^Tx\pm\mathbb{E}\left[\psi(x,\widetilde{\omega})\right]\\&\text{subject to}\quad Ax\leq b,\end{aligned} \quad (15)$$

where $P\in\mathbb{R}^{n\times n}$ is a symmetric matrix, $A\in\mathbb{R}^{\ell\times n}$, $c\in\mathbb{R}^n$ and $b\in\mathbb{R}^\ell$ are two vectors, $\mathbb{E}$ denotes the expectation of an event with respect to the random variable $\widetilde{\omega}$ defined on the probability space $(\Omega,\mathcal{F},\mathbb{P})$, with $\Omega$ being the sample space, $\mathcal{F}$ being the $\sigma$-algebra generated by subsets of $\Omega$, and $\mathbb{P}$ being a probability measure defined on $\mathcal{F}$, and $\psi(x,\omega)$ is the value function of a quadratic program given by

$$\psi(x,\omega)\triangleq\underset{y\in\mathbb{R}^m}{\text{minimize}}\left\{\frac{1}{2}y^TQy+[f(\omega)+G(\omega)x]^Ty\ \Big|\ Dy+C(\omega)x\geq\xi(\omega)\right\},$$

where $Q\in\mathbb{R}^{m\times m}$ is a symmetric matrix that may be indefinite, $D\in\mathbb{R}^{k\times m}$, $f:\Omega\to\mathbb{R}^m$ and $\xi:\Omega\to\mathbb{R}^k$ are random vector-valued functions, and $G:\Omega\to\mathbb{R}^{m\times n}$ and $C:\Omega\to\mathbb{R}^{k\times n}$ are random matrix-valued functions. The formulation (15) is a departure from traditional stochastic programming in two major ways: one is the indefinite nature of the problem that occurs in the objective functions of the two stages ($P$ not necessarily positive semidefinite in the first stage, and more interestingly, $Q$ not required to be positive semidefinite in the second stage); the other departure is the subtraction of the second-stage value function in the first-stage objective. This yields a model where the second-stage recourse function is a maximization problem that is unlike the usual case of a minimization second-stage recourse. In any event, due to its non-convexity and non-differentiability, the resulting problem (15) has not been studied in the literature of stochastic programming in this generality.

Suppose that the random variable $\omega$ is discretely distributed with $\Omega=\{\omega^1,\ldots,\omega^S\}$ for some positive integer $S$ and let $\{p_1,\ldots,p_S\}$ be the associated family of probabilities. The above stochastic quadratic program turns out to be

$$\begin{aligned}&\underset{x\in\mathbb{R}^n}{\text{minimize}}\quad \frac{1}{2}x^TPx+c^Tx\pm\sum_{s=1}^S p_s\,\psi(x,\omega^s),\\&\text{subject to}\quad Ax\leq b.\end{aligned} \quad (16)$$

By a recent result in [27, Proposition 4] about the piecewise property of the optimal objective value of the quadratic program

$$\begin{aligned}\text{val}(q,r)\triangleq\quad&\underset{y\in\mathbb{R}^m}{\text{minimum}}\quad \frac{1}{2}y^TQy+q^Ty\\&\text{subject to}\quad Dy\geq r,\end{aligned}$$

it follows that for a fixed pair $(Q, D)$ with the symmetric matrix $Q$ being copositive on the recession cone of the feasible region, i.e., $v^T Q v \geq 0$ for all $v \in \mathbb{R}^m$ satisfying $Dv \geq 0$, then $\text{val}(q, r)$ is a piecewise linear-quadratic function on the domain

$$\text{dom}(Q, D) \triangleq \left\{ (q, r) \in \mathbb{R}^{m+k} \mid -\infty < \text{val}(q, r) < +\infty \right\}.$$

Hence, under the above condition, each function $\psi(\bullet, \omega^s)$ in (16) is piecewise linear-quadratic on its domain $\mathcal{X}^s \triangleq \{x \in \mathbb{R}^n \mid -\infty < \psi(x, \omega^s) < +\infty\}$, which further yields that the objective function in (16), as the sum of piecewise linear-quadratic functions, is a piecewise linear-quadratic function [34, Exercise 10.22]. To avoid technicalities associated with extended-valued piecewise linear-quadratic functions as permitted in the reference, we state the result below by imposing a "complete recourse" assumption that allows us to directly apply Proposition 2(a).

**Proposition 15.** Suppose that $Q$ is copositive on $\mathcal{D}_\infty \triangleq \{v \in \mathbb{R}^m \mid Dv \geq 0\}$ and that $\{x \in \mathbb{R}^n \mid Ax \leq b\}$ is a subset of $\bigcap_{s=1}^{S} \mathcal{X}^s$. Then the set of d-stationary values of the two-stage stochastic quadratic program (16) is finite.

It should be noted that while there is a one-stage formulation of the problem (16) when the second-stage recourse is added to the first-stage objective in the case of finite scenarios, there is no such one-stage formulation when the second-stage recourse is subtracted in the first-stage objective. In this case, Proposition 15 adds value to the vast literature of stochastic programming with a new formulation that seems not to have been investigated.

## 5.3 Least-squares piecewise affine regression.

In this last subsection, we return to the least-squares (LS) piecewise affine regression problem that has motivated this study. Being a generalization of the classic LS linear regression problem, the least-squares piecewise affine regression problem is to find, given data points $\{x^\ell, y_\ell\}_{\ell=1}^N \in \mathbb{R}^{d+1}$, a LS estimator of a continuous piecewise affine function [9, 23, 3, 16]. Since every piecewise affine function can be written in the form (1), we may formulate this regression problem as

$$\underset{\Theta \in X}{\text{minimize}} \quad f(\Theta) \triangleq \frac{1}{2N} \sum_{\ell=1}^{N} \left( y_\ell - \left[ \max_{1 \leq i \leq k_1} \left( (a^i)^T x^\ell + \alpha_i \right) - \max_{1 \leq i \leq k_2} \left( (b^i)^T x^\ell + \beta_i \right) \right] \right)^2, \quad (17)$$

with parameter $\Theta \triangleq \left\{ (a^i, \alpha_i)_{i=1}^{k_1}, (b^i, \beta_i)_{i=1}^{k_2} \right\} \in \mathbb{R}^{(k_1+k_2)(d+1)}$ and $X$ being a polyhedral set. The above problem includes the 1-layer neural network by the rectified linear unit (ReLu) activation function [26, 14] that has the simple form $f(\theta) = \max(a^T x + \alpha, 0)$ with the parameter $\theta = (a, \alpha) \in \mathbb{R}^{d+1}$. A simple example showing the difference between the classic LS linear regression and the piecewise affine regression is demonstrated in Figure 3.

By exploiting the piecewise linear-quadratic property of the objective function, we can deduce the following result that summarizes several basic properties of the problem (17).

**Proposition 16.** The following statements hold for the problem (17).

(a) It attains a finite globally minimum value.

(b) The set of d-stationary values is finite.

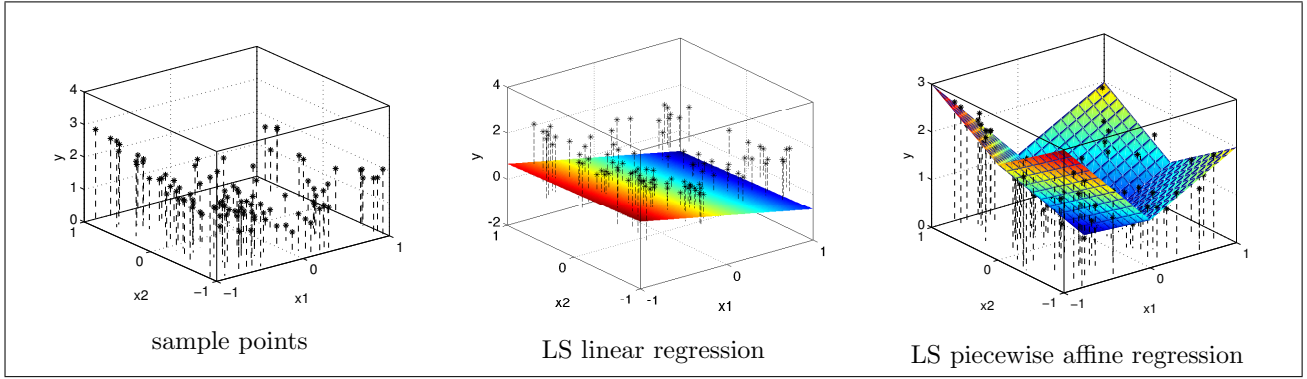

Figure 3: LS linear regression vs LS piecewise affine regression

(c) Every d-stationary point is a local minimizer.

**Proof.** Since the objective function is nonnegative, thus bounded below, statement (a) follows from the well-known Frank-Wolfe theorem for quadratic programs, extended to the case of a piecewise linear-quadratic program. Statement (b) follows from Proposition 2(a). Statement (c) holds because the objective function is the composite of a convex sum-of-squares function with a piecewise affine (max - max) function; see [7, Propositions 4.1 & 5.1] for a proof. □

In the manuscript [8], based on a convex majorization of the objective function, a modified majorization-minimization algorithm combined with a semismooth Newton method is proposed for computing a d-stationary point of a class of composite difference-max programs, which includes the problem (17) as a special case. The iterates produced by this algorithm are such that (a) the sequence of values of the majorizing functions is monotonically decreasing (but not necessarily the sequence of objective values of (17)), (b) the successive differences of the iterates converge to zero, and (c) every accumulation point of the iterates is a d-stationary solution. Thus, provided that the iterates are bounded, Proposition 14 implies that the sequence of objective values of (17) converges to a locally minimum value of this non-convex, non-differentiable optimization problem. Here, the finite number of d-stationary values is a key reason for the sequential convergence of the objective values, in spite of the lack of monotonicity of these values. For more details of this application, we refer the reader to the cited reference.

## Acknowledgments.

The work was based on research partially supported by the U.S. National Science Foundation grant IIS–1632971.